\documentclass[11pt]{amsart}
\usepackage[english]{babel}
\usepackage{egothic}
\usepackage[T1]{fontenc}

\usepackage{amsmath}
\usepackage{amsthm}
\usepackage{amssymb}
\usepackage[all]{xy}
\usepackage{mathrsfs}
\usepackage{enumerate}
\usepackage{physics}
\usepackage[notcite, final, notref]{showkeys}
\usepackage{dsfont}
\usepackage[dvips]{color}
\newtheorem{theorem}{Theorem}[section]
\newtheorem{problem}[theorem]{Problem}
\newtheorem{lemma}[theorem]{Lemma}
\newtheorem{proposition}[theorem]{Proposition}
\newtheorem{corollary}[theorem]{Corollary}
\newtheorem{definition}[theorem]{Definition}
\newtheorem{assumption}[theorem]{Assumption}

\theoremstyle{definition}
\newtheorem{remark}[theorem]{Remark}

\newtheorem{example}[theorem]{Example}

\usepackage{xcolor}

\title[]{Hörmander's $L^2$-method, $\overline{\partial}$-problem and polyanalytic function theory in one complex variable}

\author[D. Alpay]{Daniel Alpay}
\address{(DA) Schmid College of Science and Technology \\
Chapman University\\
One University Drive
Orange, California 92866\\
USA}
\email{alpay@chapman.edu}

\author[F. Colombo]{Fabrizio Colombo}
\address{(FC) Politecnico di
Milano\\Dipartimento di Matematica\\Via E. Bonardi, 9\\20133 Milano\\Italy}
\email{fabrizio.colombo@polimi.it}

\author[K. Diki]{Kamal Diki}
\address{(KD) Schmid College of Science and Technology \\
Chapman University\\
One University Drive
Orange, California 92866\\
USA}
\email{diki@chapman.edu}

\author[I. Sabadini]{Irene Sabadini}
\address{(IS) Politecnico di
Milano\\Dipartimento di Matematica\\Via E. Bonardi, 9\\20133 Milano\\Italy}
\email{irene.sabadini@polimi.it}

\author[D. C. Struppa]{Daniele C. Struppa}
\address{(DCS) Donald Bren Presidential Chair in Mathematics \\
Chapman University\\
One University Drive
Orange, California 92866\\
USA}
\email{struppa@chapman.edu}

\begin{document}
\maketitle
\tableofcontents
\begin{abstract}
In this paper we consider the classical $\overline{\partial}$-problem in the case of one complex variable both for analytic and polyanalytic data. We apply the decomposition property of polyanalytic functions in order to construct particular solutions of this problem and obtain new Hörmander type estimates using suitable powers of the Cauchy-Riemann operator. We  also compute particular solutions of the $\overline{\partial}$-problem for specific polyanalytic data such as the It{\^o} complex Hermite polynomials and polyanalytic Fock kernels.

\end{abstract}


\section{Introduction}
\setcounter{equation}{0}

An interesting generalization of the classical theory of holomorphic functions of a complex variable is the theory of polyanalytic functions, defined as nullsolutions of higher order powers of the Cauchy-Riemann operator. These functions play an important role since both the complex variable $z$ and its conjugate $\overline{z}$ are involved in their construction, allowing  the consideration of polynomials  $P(z,\overline{z})$ of two connected variables $z$, $\bar z$.  These functions
were first introduced in 1908 by Kolossoff in \cite{kolossoff1908problemes} to study elasticity problems.
 A complete introduction to polyanalytic functions and their basic properties can be found in \cite{Balk1991, balk1997polyanalytic}. This class of functions was studied by various authors from different perspectives, see \cite{ abreu2010sampling, abreu2012super, abreu2014function, agranovsky2011characterization, AhernBruna, begehr2002orthogonal, vasilevski1999structure}
 and the references therein.
In quantum mechanics polyanalytic functions are also relevant for the study of the Landau levels associated to
Schrödinger operators, see \cite{abreu2014function, askour2000explicit}. These functions were used also in \cite{abreu2010sampling} to study sampling and interpolation problems
on polyanalytic Fock spaces using time frequency analysis techniques such as short-time Fourier transform or
Gabor transforms.

\vskip0.2cm

In this paper we use ideas and techniques from the theory of polyanalytic functions to construct particular and general solutions for the classical $\overline{\partial}$-problem, as well as for the $\overline{\partial}$-problem in which the datum is a polyanalytic function. The main tools that we will use are the well known $L^2$ estimates that H\"ormander developed, see for example \cite{hormander1966l2} and \cite{hormander1967generators, hormander2007notions}.

\vskip0.2cm

The structure of the paper is as follows: in Section 2 we review some well-known facts on the $\overline{\partial}$-problem and recall H\"ormander's classical $L^2$-estimates for the case of one complex variable. We recall also the definition of polyanalytic functions and present some of their main properties, introducing as well a notion of Fock space in the polyanalytic setting. In Section 3 we treat the case of the $\overline{\partial}$-problem for a given analytic datum. We give an equivalent representation of the spaces of entire functions with exponential growth and provide some interesting examples which appear naturally from the computations. Section 4 is devoted to the case where the given datum is a polyanalytic function of order $n>1$. Here also we construct a particular solution of the $\overline{\partial}$-problem using a suitable differential operator of finite order involving powers of the conjugate multiplication operator and the classical Cauchy-Riemann operator. We obtain new $L^2$-estimates in this framework involving a Sobolev-type norm. Some particular examples of polyanalytic data including complex Hermite polynomials and polyanalytic Fock spaces are studied as applications of this approach.

\section{Preliminary notations and results}
\setcounter{equation}{0}
We begin by discussing a modified version of the Hörmander's result in the one complex variable setting as presented in \cite[Theorem 1]{berenstein1979new}(see also \cite[Theorem 1.4]{berenstein1993complex} for the case of several complex variables). Let $\Omega\subset \mathbb{C}$ be an open set; we denote by $W_p(\Omega)$ the space of all measurable functions $f:\Omega \longrightarrow \mathbb{C}$ for which there exists a constant $C>0$ such that
%\begin{theorem}\label{HörmanderC}
%Let $\varphi$ be a plurisubharmonic function on $\mathbb{C}^n$, and let $f$ be a $(0,1)$-form such that $\overline{\partial}f=0$. If, moreover,

%\begin{equation}
%\displaystyle \int_{\mathbb{C}^n}|f(z)|^2e^{-2\varphi(z)}d\lambda(z)<+\infty,
%\end{equation}
%with $d\lambda$ being the Lebesgue measure on $\mathbb{C}^n,$ then there exists a function $u$ on $\mathbb{C}^n$ such that
%\begin{equation}
%\displaystyle \int_{\mathbb{C}^n}\frac{|u(z)|^2}{(1+|z|^2)^2}e^{-2\varphi(z)}d\lambda(z)\leq \int_{\mathbb{C}^n}|f(z)|^2e^{-2\varphi(z)}d\lambda(z),
%\end{equation}
%and which is a solution of the equation
%\begin{equation}
%\overline{\partial}u=f.
%\end{equation}
%\end{theorem}

\begin{equation}
\displaystyle \int_{\Omega}|f(z)|^2e^{-Cp(z)}d\lambda(z)<+\infty,
\end{equation}
where $d\lambda(z=x+iy)=dxdy$ is the classical Lebesgue measure on $\mathbb{C}$ and $p$ is a subharmonic function. Throughout the entirety of this paper we will always assume that $p$ satisfies the following:
\begin{assumption}\label{as12}
\ \
\begin{enumerate}
\item $p(z)\geq 0$ and $\log(1+|z|^2)=O(p(z));$
\item there exist constants $C, D>0$ such that $|z -w|\leq 1$ implies that $p(z)\leq Cp(w)+D$.
\end{enumerate}
\end{assumption}
Theorem 1 of \cite{berenstein1979new} can be stated as follows:
\begin{theorem}\label{HormanderBT}
Let $\Omega$ be an open subset of $\mathbb{C}$, let $p$ be a subharmonic function in $\Omega$, and let $f$ be a function in $W_{p}(\Omega)$ such that
\begin{equation}
\displaystyle \int_{\Omega}|f(z)|^2e^{-p(z)}d\lambda(z)=M(f)<+\infty.
\end{equation}
 Then, there exists a solution $u$ to equation
 \begin{equation}
\overline{\partial}u:=\frac{\partial}{\partial \overline{z}}u=f
\end{equation}
satisfying
\begin{equation}
\displaystyle \int_{\Omega}\frac{|u(z)|^2}{(1+|z|^2)^2}e^{-p(z)}d\lambda(z)\leq \frac{M(f)}{2}.
\end{equation}
Moreover, if the function $f$ belongs to $\mathcal{C}^{\infty}(\Omega)$ then the solution $u$ belongs to $\mathcal{C}^{\infty}(\Omega)$.
\end{theorem}
We are also interested in the space 
\begin{equation}\label{Apspace}
\displaystyle A_p= A_p(\mathbb{C})=\left\lbrace{f\in H(\mathbb{C});\quad \exists A,B>0 \quad |f(z)|\leq A\exp(Bp(z))}\right\rbrace
\end{equation}
of entire functions whose growth is controlled by $p$. Note that property (1) implies that all the polynomials belong to the space $A_p$. Furthermore property (2) implies that the space $A_p$ is stable by differentiation and finally
\begin{equation}\label{HormRes}
A_p(\mathbb{C})=H(\mathbb{C})\cap W_p(\mathbb{C}).
\end{equation}
Let us turn our attention to the theory of polyanalytic functions of one complex variable, following \cite{Balk1991}. A complex-valued function $f:\Omega\subset \mathbb{C}\longrightarrow \mathbb{C}$ of class $\mathcal{C}^n$ on a domain $\Omega$  is called a polyanalytic function of order $n$ if it belongs to the kernel of the $n-$th power, $n\geq 1,$ of the classical Cauchy-Riemann operator $\displaystyle \frac{\partial}{\partial \overline{z}}$, that is if $$\displaystyle \frac{\partial^n}{\partial \overline{z}^n}f(z)=0, \quad \forall  z\in\Omega.$$
\noindent 
We denote the space of polyanalytic functions of order $n$ by $H_n(\Omega)$. When $n=1$ we get the space of analytic functions and we write $H(\Omega)$ instead of $H_1(\Omega)$.
\smallskip
An interesting fact regarding these functions is that any polyanalytic function of order $n$ can be decomposed in terms of $n$ analytic functions so that we have a decomposition of the following form
\begin{equation}\label{polydecomposition}
f(z)=\displaystyle \sum_{k=0}^{n-1}\overline{z}^kf_k(z),
\end{equation}
for which all $f_k$ are analytic functions on $\Omega$. Expanding each analytic component using the series expansion at the origin leads to an expression of the form

\begin{equation}\label{exp1}
f(z)=\displaystyle \sum_{k=0}^{n-1}\sum_{j=0}^{\infty}\overline{z}^kz^ja_{k,j},
\end{equation}
where $(a_{k,j})$ are complex coefficients.

Another space which will be useful in the sequel is the Fock space, that we introduce in the more general context of polyanalytic functions. For $n=1,2,...$ the polyanalytic Fock spaces of order $n$ are defined as
$$\mathcal{F}_n(\mathbb{C}):=\left\lbrace g\in H_n(\mathbb{C}): \quad \frac{1}{\pi}\int_{\mathbb{C}}|g(z)|^2e^{-|z|^2}d\lambda(z)<\infty \right\rbrace.$$
For $n=1$ we obviously get the classical Fock space.\\
The reproducing kernel associated to the space $\mathcal{F}_n(\mathbb{C})$ is given by
\begin{equation}\label{Kn}
F_n(z,w)=e^{z\overline{w}}\displaystyle \sum_{k=0}^{n-1}\frac{(-1)^k}{k!}{n \choose k+1}|z-w|^{2k},
\end{equation}
for every $z,w\in\mathbb{C}.$
\section{Hörmander's $L^2$-method: polyanalytic trick}
\setcounter{equation}{0}
In this section we will apply H\"ormander's $L^2-$estimates from Theorem \ref{HormanderBT}, combined with the approach proposed by Berenstein and Taylor in \cite{berenstein1979new},  to the case of polyanalytic function theory.
\subsection{Characterization of $\overline{\partial}$-problem solutions: analytic case}

 Our goal here is to use the theory of polyanalytic functions to give an explicit characterization of the space of solutions to the $\overline{\partial}$ problem. We first introduce the following space:

\begin{definition}\label{Hspace}
Let $f\in W_p(\Omega)$ be such that $$\int_{\Omega}|f(z)|^2e^{-p(z)}d\lambda(z)=M(f)<+\infty.$$  We define the set
$$\mathcal{H}_{f,p}(\Omega) =\left\lbrace{g\in H(\Omega); \quad \int_{\Omega}\frac{|g(z)|^2}{(1+|z|^2)^2}e^{-p(z)}d\lambda(z)\leq 3M(f)}\right\rbrace.$$
\end{definition}

We now can prove the first result of this paper.

\begin{theorem}\label{H-poly}
Let $\Omega$ be an open set of $\mathbb{C}$,  let $p$ be subharmonic in $\Omega$, and let $f \in W_{p}(\Omega)$ be such that $\overline{\partial}f=0$ and \begin{equation}
\displaystyle \int_{\Omega}|f(z)|^2e^{-p(z)}d\lambda(z)<+\infty.
\end{equation}
Then there exists a polyanalytic function  $u$ of order $2$ which is a solution of the problem $$\overline{\partial}u=f,$$ and which can be expressed as follows
\begin{equation}\label{uexpress}
u(z,\overline{z})=\overline{z}f(z)+u_0(z), \quad
\end{equation}
where $u_0$ belongs to $\mathcal{H}_{f,p}(\Omega)$.
\end{theorem}
\begin{proof}
Since $f\in W_p(\Omega)$, Theorem \ref{HormanderBT} implies the existence of a solution $u$  to the equation
\begin{equation}
\overline{\partial}u=f,
\end{equation}
which satisfies the estimate
\begin{equation}\label{EST1}
\displaystyle \int_{\Omega}\frac{|u(z,\bar z)|^2}{(1+|z|^2)^2}e^{-p(z)}d\lambda(z)\leq \frac{1}{2}\int_{\Omega}|f(z)|^2e^{-p(z)}d\lambda(z)=\frac{M(f)}{2}.
\end{equation}

Since $f\in H(\Omega)$, it is in particular of class $\mathcal{C}^{\infty}$, and therefore Theorem \ref{HormanderBT} implies that $u$ is of class $\mathcal{C}^{\infty}$ as well. Thus, we can apply the $\overline{\partial}$-operator and get
$$\overline{\partial}^2u(z)=\overline{\partial}f(z)=0,\quad \forall z\in \Omega.$$
Hence, $u$ is polyanalytic of order $2$ everywhere on $\Omega$. Then, by the poly-decomposition (\ref{polydecomposition}), $u$ can be expressed in terms of two holomorphic functions $u_0,u_1\in \mathcal{H}(\Omega)$ so that
\begin{equation}\label{udecom}
u(z,\overline{z})=u_0(z)+\overline{z}u_1(z), \quad  z\in\Omega.
\end{equation}
We now apply the $\overline{\partial}$-operator to expression \eqref{udecom}, and using the fact that both $u_0$ and $u_1$ are in the kernel of $\overline{\partial}$ we obtain
$$\overline{\partial}u(z,\overline{z})=\overline{\partial}u_0(z)+\overline{z}\overline{\partial}u_1(z)+u_1(z)=u_1(z).$$

Moreover, since $\overline{\partial}u=f$ it is clear that $u_1=f$ everywhere on $\Omega$. Hence, we obtain
$$u(z,\overline{z})=\overline{z}f(z)+u_0(z), \quad \forall z\in\Omega.$$
This implies in particular that $$u_0(z)=u(z,\overline{z})-\overline{z}f(z).$$

This implies that $$|u_0(z)|^2\leq 2(|u(z,\overline{z})|^2+|z|^2|f(z)|^2), \quad \forall z\in\Omega.$$
We continue the computations and get
\begin{equation}\label{EST2}
\displaystyle \int_{\Omega}\frac{|u_0(z)|^2}{(1+|z|^2)^2}e^{-p(z)}d\lambda(z)\leq 2\left(\int_{\Omega}\frac{|u(z,\overline{z})|^2}{(1+|z|^2)^2}e^{-p(z)}d\lambda(z)+ \int_{\Omega}\frac{ |z|^2|f(z)|^2}{(1+|z|^2)^2}e^{-p(z)}d\lambda(z)\right).
\end{equation}

Since $$|z|^2\leq 1+|z|^2\leq (1+|z|^2)^2, \quad \forall z\in\Omega,$$
we have $$\displaystyle \int_{\Omega}\frac{ |z|^2|f(z)|^2}{(1+|z|^2)^2}e^{-p(z)}d\lambda(z)\leq \int_{\Omega}|f(z)|^2e^{-p(z)}d\lambda(z)=M(f).$$

Then, we use estimate \eqref{EST1} and insert it in \eqref{EST2} to deduce that
\begin{equation}
\displaystyle \int_{\Omega}\frac{|u_0(z)|^2}{(1+|z|^2)^2}e^{-p(z)}d\lambda(z)\leq 3 \int_{\Omega} |f(z)|^2e^{-p(z)}d\lambda(z)=3M(f)<+\infty.
\end{equation}
Finally we conclude that $u$ can be expressed as in \eqref{uexpress} with $u_0$ belonging to $\mathcal{H}_{f,p}(\Omega)$.
\end{proof}
%Conversely, we can prove the following
\begin{proposition}
Let $\Omega$ be an open subset of $\mathbb{C}$,  let $p$ be subharmonic in $\Omega$ and let $f \in W_{p}(\Omega)$ be such that $\overline{\partial}f=0$ and \begin{equation}
\displaystyle \int_{\Omega}|f(z)|^2e^{-p(z)}d\lambda(z)=M(f)<+\infty.
\end{equation}

Assume that
\begin{equation}
u(z,\overline{z})=\overline{z}f(z)+u_0(z),
\end{equation}
for every $z\in\Omega$ with $u_0\in \mathcal{H}_{f,p}(\Omega)$. Then:
\begin{enumerate}
\item[i)] $u$ solves the problem $\overline{\partial}u=f;$
\item[ii)] $u$ satisfies 
 $$\displaystyle \int_{\Omega}\frac{|u(z,\overline{z})|^2}{(1+|z|^2)^2}e^{-p(z)}d\lambda(z)\leq 8M(f)<+\infty.$$
\end{enumerate}
\end{proposition}
\begin{proof}
(i) We note that $u_0$ is holomorphic on $\Omega$, so it is clear that $u$ is a particular solution in this case. In fact, we have $$\overline{\partial}u=\overline{\partial}(\overline{z}f)+\overline{\partial}(u_0)=f.$$
\\
(ii) Since $u_0\in \mathcal{H}_{f,p}(\Omega)$ and $$|u(z)|^2\leq 2(|u_0(z,\overline{z})|^2+|z|^2|f(z)|^2), \quad \forall z\in\Omega,$$
we have
$$\int_{\Omega}\frac{|u_0(z)|^2}{(1+|z|^2)^2}e^{-p(z)}d\lambda(z)\leq 3M(f),$$
as well as
$$\int_{\Omega}\frac{|z|^2|f(z)|^2}{(1+|z|^2)^2}e^{-p(z)}d\lambda(z)\leq M(f).$$ It is then easy to check that
$$\displaystyle \int_{\Omega}\frac{|u(z,\overline{z})|^2}{(1+|z|^2)^2}e^{-p(z)}d\lambda(z)\leq 8M(f)<+\infty.$$

\end{proof}
As a consequence of the previous results we have the following:
\begin{corollary}
Let $p$ be a subharmonic function in $\mathbb{C}$.
If $f$ belongs to $A_p(\mathbb{C})$, and
$$\displaystyle \int_{\mathbb{C}}|f(z)|^2e^{-p(z)}d\lambda(z)<+\infty,$$
  then every solution $u$ of the equation $\overline{\partial}u=f$ can be expressed as $$u(z,\overline{z})=\overline{z}f(z)+u_0(z),$$
with $u_0\in \mathcal{H}_{f,p}(\mathbb{C})$.
\end{corollary}
\begin{proof}
We apply Theorem \ref{H-poly} to the case of entire functions by taking $\Omega=\mathbb{C}$ and using \eqref{HormRes}.
\end{proof}
%\begin{remark} Let $f\in H(\Omega)$,  we can consider also a general space with estimate given by
 %$$\mathbf{H}_{f}=\left\lbrace{g\in H(\Omega); \quad |g(z)| \leq 2(1+|z|^2)|f(z)|}\right\rbrace.$$

%It is easy to see that $$\mathbf{H}_{f}\subset \mathcal{H}_{f,p},$$
%where $p$ is a subharmonic function on $\mathbb{C}$.
%\end{remark}
We now give a first example based on the classical Fock space $\mathcal{F}(\mathbb{C})$. We take $p(z)=|z|^2$ and $\Omega=\mathbb{C}$.
\begin{example}
For every fixed parameter $w\in \mathbb{C}$, we set 
\begin{equation}\label{Kn1}
F(z,w)=F_w(z)=e^{z\overline{w}}, \quad \forall z\in \mathbb{C}.
\end{equation}
Then $u$ is a solution of the $\overline{\partial}-$problem with datum $F_w$ if and only if
$$u(z,\overline{
z})=\overline{z}e^{z\overline{w}}+u_0(z),\quad \forall z\in\mathbb{C},$$
with $u_0\in\mathcal{H}_{F_w}(\mathbb{C}).$
In fact, this a direct consequence of Theorem \ref{H-poly} by taking $\Omega=\mathbb{C}$ and $f(z)=F_w(z)$.
\end{example}
\begin{remark}
Note that in the previous example we have
$$M(F_w)=||F_w||_{\mathcal{F}(\mathbb{C})}^2=e^{|w|^2};\quad \forall w\in \mathbb{C}.$$
\end{remark}
\subsection{Reproducing kernel Hilbert space associated to $\overline{\partial}$-problem: Gaussian case}
In this subsection, we are interested in the case $p(z)=|z|^2$ and $\Omega=\mathbb{C}$. 
%Then we consider a datum $f$ which belongs to the classical Fock space $\mathcal{F}(\mathbb{C})$, i.e. $\overline{\partial }f=0$ on $\mathbb{C}$ and such that
%$$M(f):=\displaystyle \int_{\mathbb{C}}|f(z)|^2e^{-|z|^2}d\lambda(z)<+\infty.$$
%\begin{example}
%If $\varphi(z)=\frac{|z|^2}{2}$, the previous result presents solutions of the $\overline{\partial}$-problem for the classical Fock-Bargmann space.
%\end{example}
We begin by introducing the space $\mathcal{H}_p$, consisting of entire functions and containing the functions $u_0$ as in \eqref{uexpress}. 
\begin{definition}\label{Hpdefi}
Let $p$ be a subharmonic function on $\mathbb{C}$. We define the space
$$\mathcal{H}_p:=\left\lbrace{g\in H(\mathbb{C});\quad ||g||^2_{\mathcal{H}}:=\frac{1}{\pi}\int_{\mathbb{C}}\frac{|g(z)|^2}{(1+|z|^2)^2}e^{-p(z)}d\lambda(z)<+\infty}\right\rbrace,$$
with inner product defined by \begin{equation}
\displaystyle \langle u,v \rangle_{\mathcal{H}_p}:=\frac{1}{\pi}\int_{\mathbb{C}}\frac{u(z)\overline{v(z)}}{(1+|z|^2)^2}e^{-p(z)}d\lambda(z).
\end{equation}
\end{definition}
\begin{remark}
The space $\mathcal{H}_p$ introduced above corresponds to the space of functions considered by Hörmander in \cite[Theorem 2.5.3]{hormander1966l2}.
\end{remark}
\begin{definition}[Hörmander-Fock space] \label{HFS}
When $p(z)=|z|^2$, the space $\mathcal{H}_{p}$ will be denoted simply by $\mathcal{H}$ and will be called the Hörmander-Fock type space.
\end{definition}
\begin{proposition}
The Fock space $\mathcal{F}(\mathbb{C})$ is included in the Hörmander-Fock space $\mathcal{H}$. In other words, the injection $$\mathcal{F}(\mathbb{C})\hookrightarrow \mathcal{H},$$ is continuous, so that for every $g\in \mathcal{F}(\mathbb{C})$, it holds \begin{equation}
||g||_{\mathcal{H}}\leq ||g||_{\mathcal{F}(\mathbb{C})}.
\end{equation}
\end{proposition}
\begin{proof}
The assertion follows from the fact that fhe following estimate holds for any $g\in\mathcal{F}(\mathbb{C})$  $$||g||^2_{\mathcal{H}}:=\int_{\mathbb{C}}\frac{|g(z)|^2}{(1+|z|^2)^2}e^{-|z|^2}d\lambda(z)\leq  \int_{\mathbb{C}}|g(z)|^2e^{-|z|^2}d\lambda(z)=||g||_{\mathcal{F}(\mathbb{C})}<+\infty.$$
\end{proof}

The spaces $\mathcal{F}(\mathbb{C})$ and $\mathcal{H}(\mathbb{C})$ are different as the following example shows.
\begin{example}
Since the creation operator $M_z(f)(z)=zf(z)$ is unbounded on the Fock space $\mathcal{F}(\mathbb{C})$, we can choose a function $f\in\mathcal{F}(\mathbb{C})$ such that the function $g=M_z(f)=zf\notin \mathcal{F}(\mathbb{C})$. It is clear that we have the following estimate
$$\displaystyle \int_\mathbb{C}\frac{|g(z)|^2}{(1+|z|^2)^2}e^{-|z|^2}d\lambda(z)\leq ||f||_{\mathcal{F}(\mathbb{C})}<+\infty.$$
Thus $g$ belongs to the Hörmander-Fock space $\mathcal{H}$ but does not belong to the Fock space $\mathcal{F}(\mathbb{C})$. Hence, the Fock space $\mathcal{F}(\mathbb{C})$ is strictly included in the space $\mathcal{H}$.
\end{example}
The Hörmander-Fock type space $\mathcal H$ is of independent importance, see \cite{HFspace}. In particular, we have:
\begin{proposition}
For every $n=0,1,...$ the norm of the monomials $z^n$ with respect to the Hilbert space $\mathcal{H}$ is given by the following moment sequence
\begin{equation}
\displaystyle \eta_n=||z^n||_{\mathcal{H}}^2=\int_0^{\infty}\frac{t^n}{(1+t)^2}e^{-t}dt\leq n!.
\end{equation}
Moreover, we have the orhogonality condition
\begin{equation}
\langle z^n,z^m \rangle_{\mathcal{H}}=0, \quad \forall n\neq m\in \mathbb{N}.
\end{equation}
\end{proposition}
The moment sequence $\eta_n$ gives rise to a number of results, which are the main objective of \cite{HFspace}.

\subsection{An equivalent representation of the spaces $A_p(\mathbb{C})$}
In this subsection we show that the space $$\displaystyle A_p= A_p(\mathbb{C})=\left\lbrace{f\in H(\mathbb{C});\quad \exists A,B>0 \quad |f(z)|\leq A\exp(Bp(z))}\right\rbrace,$$

coincides with the space $\mathfrak{H}_p$ consisting of all entire functions $f:\mathbb{C}\longrightarrow \mathbb{C}$ such that there exists a constant $C>0$ for which we have
\begin{equation}
\displaystyle \int_\mathbb{C} \frac{|f(z)|^2}{(1+|z|^2)^2}e^{-Cp(z)}d\lambda(z)<+\infty.
\end{equation}

\begin{theorem}
An entire function $f:\mathbb{C}\longrightarrow \mathbb{C}$ belongs to the space $A_p(\mathbb{C})$ if and only if $f$ belongs to $\mathfrak{H}_p(\mathbb{C}).$  In other words, as sets, we have
$$A_p(\mathbb{C})=\mathfrak{H}_p(\mathbb{C}).$$
\end{theorem}
\begin{proof}
If we assume that $f\in A_p(\mathbb{C})$ we have the estimate
\begin{equation}
\displaystyle \int_{\mathbb{C}} \frac{|f(z)|^2}{(1+|z|^2)^2}e^{-Cp(z)}d\lambda(z)\leq \int_{\mathbb{C}} |f(z)|^2e^{-Cp(z)}d\lambda(z)<+\infty.
\end{equation}
So, in particular,
 $$A_p(\mathbb{C})\subset \mathfrak{H}_p(\mathbb{C}).$$

 For the converse, let us take $f\in \mathfrak{H}_p(\mathbb{C})$ and show that $f$ belongs to $A_p(\mathbb{C})$. Setting

 \begin{equation}\label{qformula}
  q(z)=Cp(z)+2\log(1+|z|^2),
\end{equation}

 we observe that since $f$ belongs to $\mathfrak{H}_p(\mathbb{C})$ we have
 \begin{equation}
 \displaystyle \int_{\mathbb{C}} |f(z)|^2e^{-q(z)}d\lambda(z)=\int_{\mathbb{C}} \frac{|f(z)|^2}{(1+|z|^2)^2}e^{-Cp(z)}d\lambda(z)<+\infty.
 \end{equation}

 Thus we easily deduce that  $f\in A_q(\mathbb{C})$, so that there exist $A, B>0$ such that
 $$|f(z)| \leq A\exp(Bq(z)), \quad \forall z\in \mathbb{C}.$$
However, since $\log(1+|z|^2)=O(p(z))$, there exists $M>0$ such that we have $$\log(1+|z|^2)\leq Mp(z), \quad \forall z\in \mathbb{C}.$$
So, it is clear from the definition of the polynomial $q(z)$ given by formula \eqref{qformula} that
$$q(z)\leq (C+2M)p(z), \quad \forall z\in \mathbb{C}.$$
Setting $B'=B(C+2M)$ we obtain the estimate below
\begin{align*}
\displaystyle |f(z)|&\leq A \exp\left(B(Cp(z)+2\log(1+|z|^2)) \right) \\
&\leq A\exp(B'p(z)). \\
\end{align*}

So, we deduce that we always can find two constants $A', B'>0$ such that
 $$|f(z)| \leq A'\exp(B'p(z)), \quad \forall z\in \mathbb{C}.$$

 In particular this shows that $f$ belongs to $A_p(\mathbb{C})$, so we have

$$\mathfrak{H}_p(\mathbb{C})\subset A_p(\mathbb{C}).$$
 This concludes the proof.

%$$Cp(z)\leq q(z)$$

%However, $p(z)=O(\log(1+|z|^2))$ and $p(z)\geq 0$, so, there exists $M>0$ such that
%$$p(z)\leq M \log(1+|z|^2).$$
%Hence, in particular we have $q(z)\geq 0$ and
%$$q(z)\leq (M+2)\log(1+|z|^2)\leq (M+2)|z|^2; \quad \forall z.$$
%Then, by setting $K=2+M$ we have
%$$e^{-\frac{1}{K}q(z)}\geq e^{-|z|^2}$$
\end{proof}
\begin{remark}
The space $\mathcal{H}_p(\mathbb{C})$ considered in Definition \ref{Hpdefi} is contained in the space $A_p(\mathbb{C})$ but we do not have $\mathcal{H}_p= A_p$ since in the definition of $\mathcal{H}_p$ the constant is fixed to be $C=1$.
\end{remark}

\section{$\overline{\partial}$-problem and Hörmander's estimate for a polyanalytic datum}
\setcounter{equation}{0}
In this section we want to characterize the solutions of $\overline{\partial}u=f$ when the datum $f$ is a polyanalytic function of arbitrary order.  We begin with a simple example:
\begin{example}
Let $f:\mathbb{C}\longrightarrow \mathbb{C}$ be a polyanalytic function of order $2$. Setting $g(z)=\overline{z}f(z)$, an easy computation shows that
$\overline{\partial} g=f+\overline{z} \overline{\partial}(f)$; thus, unlike in the analytic case considered in Section 3, the function $g$ is not anymore a particular solution of the equation
$
\displaystyle \overline{\partial} u=f$ since the term $\overline{\partial}(f)$ is not zero.

\end{example}

We state now the more general problem that we study in this part of the paper :
\begin{problem}\label{pb-p}
Let $\Omega\subseteq \mathbb{C}$ be a domain and let $f:\Omega \longrightarrow \mathbb{C}$ be a polynalytic function of order $n>1$, i.e. $\overline{\partial}^n f(z)=0,$ for every $z\in \Omega.$
Consider the following $\overline{\partial}-$problem
\begin{equation}\label{pb-poly}
\displaystyle \overline{\partial} u(z,\overline{z})=f(z,\overline{z}).
\end{equation}
 Can we characterize particular and general solutions of the problem \eqref{pb-poly} in terms of the given polyanalytic datum $f$ and of the powers $\overline{\partial}^k(f)$ where $k=1,\cdots, n-1$ ?
\end{problem}

\subsection{The case of a bi-analytic datum}
We first start by treating the case where the given datum $f$ in Problem \ref{pb-p} is a polyanalytic function of order $2$. We first introduce the following function space:
\[
\begin{split}
\mathcal{H}^{2}_{f,p}(\Omega)&=\left\{g\in H(\Omega); \exists C>0:  \right.
\\
&\left.\int_{\Omega}\frac{|g(z)|^2}{(1+|z|^2)^4}e^{-p(z)}d\lambda(z)\leq C\displaystyle \int_{\Omega}\left(|f(z)|^2+ \left|\overline{\partial}f(z)\right|^2\right)e^{-p(z)}d\lambda(z)<+\infty\right\}.
\end{split}
\]
\begin{theorem}\label{H-polydatum}
Let $\Omega$ be an open set of $\mathbb{C}$,  let $p$ be a subharmonic function in $\Omega$ and let $f:\Omega\longrightarrow \mathbb{C}$ be a polyanalytic function of order $2$ such that \begin{equation}
\displaystyle \int_{\Omega}|f(z)|^2e^{-p(z)}d\lambda(z)<+\infty, \quad \int_{\Omega}\left|\overline{\partial}f(z)\right|^2e^{-p(z)}d\lambda(z)<+\infty.
\end{equation}
Then there exists a polyanalytic function $u$  of order $3$ which is a solution of the problem $\overline{\partial}u=f$ and which can be expressed as 
\begin{equation}\label{upoly2}
u(z,\overline{z})=-\overline{z}^2\overline{\partial}f/2+\overline{z}f+u_0, \quad
\end{equation}
where $u_0\in\mathcal{H}^{2}_{f,p}(\Omega)$.
\end{theorem}

\begin{proof}

Theorem \ref{HormanderBT} asserts the existence of a solution $u$ solving the equation $\overline{\partial}u=f$. Applying twice the $\overline{\partial}$ operator and using the fact that $f$ is polyanalytic of order $2$ we obtain
$$\overline{\partial}^3 u=\overline{\partial}^2f=0.$$
Thus, the solution $u$ will be polyanalytic of order $3$ and therefore
\begin{equation}
u(z,\overline{z})=u_0+\overline{z}u_1+\overline{z}^2u_2;
\end{equation}
where $u_0,u_1$ and $u_2$ all belong to $H(\Omega)$.
Immediate computations show that
$$\displaystyle \overline{\partial}u=u_1+2\overline{z}u_2=f;$$
and $$\overline{\partial}^2u=2u_2=\overline{\partial}f.$$

Thus, we deduce that
$$u_1=f-\overline{z}\overline{\partial}f, \quad u_2=\frac 12 \overline{\partial}f,$$ which are both holomorphic functions and therefore
\begin{equation}\label{upoly2}
u(z,\overline{z})=\frac{\overline{z}^2}{2}\overline{\partial}f(z)+\overline{z}(f(z)-\overline{z}\overline{\partial}f(z))+u_0(z), \quad \forall z\in \Omega,
\end{equation}
so that $$u(z,\overline{z})=-\frac{\overline{z}^2}{2}\overline{\partial}f+\overline{z}f(z)+u_0(z),\quad \forall z\in \Omega.$$
In order to obtain an estimate for the function $u_0$ we use the previous formula and we get
\begin{align*}
\displaystyle |u_0(z)|^2&\leq 2\left(|u(z,\overline{z})|^2+\left|\frac{\overline{z}^2}{2}\overline{\partial}(f)-\overline{z}f\right |^2\right) \\
&\leq 2|u(z,\overline{z})|^2+|z|^4|\overline{\partial}f|^2+4|z|^2|f(z)|^2 \\
\end{align*}
Therefore, we deduce
\begin{align*}
\displaystyle \int _\Omega \frac{|u_0(z)|^2}{(1+|z|^2)^4} e^{-p(z)}d\lambda(z)&\leq \int_\Omega \left( 2\frac{|u(z,\overline{z})|^2}{(1+|z|^2)^4}+\frac{|z|^4|\overline{\partial}f|^2}{(1+|z|^2)^4}+4 \frac{|z|^2|f(z)|^2}{(1+|z|^2)^4}\right) e^{-p(z)}d\lambda(z) \\
\end{align*}
Further  developing the computations using the $L^2$-estimate of $u$ we see that
\begin{align*}
\displaystyle \int _\Omega \frac{|u_0(z)|^2}{(1+|z|^2)^4} e^{-p(z)}d\lambda(z)&\leq 5 \left( \displaystyle \int_{\Omega}|f(z)|^2e^{-p(z)}d\lambda(z)+ \int_\Omega\left|\overline{\partial}f(z)\right|^2 e^{-p(z)}d\lambda(z)\right) <+\infty .\\
\end{align*}
So, we conclude that $u_0(z)$ belongs to the space $\mathcal{H}^2_{f,p}(\Omega)$, which concludes the proof.
\end{proof}
\begin{example}
Let $F_2(z,w)$ be the reproducing kernel of the polyanalytic Fock space of order $2$, see (\ref{Kn}). For every fixed parameter $w\in \mathbb{C}$ we set
\begin{equation}\label{Kn2}
g_w(z)=F_2(z,w)=F_{2,w}(z)=e^{z\overline{w}}\displaystyle \left(2-|z-w|^{2}\right), \quad \forall z\in \mathbb{C}.
\end{equation}

We can show that a particular solution of the problem $\overline{\partial}u=g_w$ is given by
\begin{equation} \label{upF2}
u(z)=\overline{z}F_2(z,w)+\frac{\overline{z}^2}{2}(z-w)e^{z\overline{w}}=\overline{z}F_2(z,w)+\frac{\overline{z}^2}{2}(z-w)F_1(z,w),\quad \forall z,w \in\mathbb{C}.
\end{equation}

Indeed, we observe that $g_w$ is polyanalytic of order $2$ with respect to the variable $z$, so that we have $\overline{\partial}^2g_w(z)=0$. By Theorem \ref{H-polydatum} we know that a particular solution is given by
$$
u(z,\overline{z})=-\frac{\overline{z}^2}{2}\overline{\partial} g_w(z)+\overline{z}g_w(z), \quad \forall z\in \mathbb{C}.
$$

It is easy to check that $$\overline{\partial}g_w(z)=-(z-w)e^{z\overline{w}}=-(z-w)F_1(z,w),$$
and this leads to the expression given by formula \eqref{upF2}.
\end{example}

\subsection{General polyanalytic datum}
In this subsection we start from a general polyanalytic datum $f$. We show how the solutions look for the cases $n=1, 2, 3, 4,$ and we then offer the general formula. These first few cases can be easily checked, so we omit the details.
\begin{example}\label{Expn4}
Here are the solutions to Problem \ref{pb-p} for $n=1, 2, 3, 4.$
\begin{itemize}
\item[i)] If $\overline{\partial}f=0,$ we have
$$u(z,\overline{z})=u_0(z)+\overline{z}f(z);$$
\item[ii)] If $\overline{\partial}^2f=0,$ we have
$$u(z,\overline{z})=u_0(z)+\overline{z}f(z)-\frac{\overline{z}^2}{2}\overline{\partial}(f);$$
\item[iii)] If $\overline{\partial}^3f=0,$ we have
$$u(z,\overline{z})=u_0(z)+\overline{z}f(z)-\frac{\overline{z}^2}{2}\overline{\partial}(f)+\frac{\overline{z}^3}{6}\overline{\partial}^2(f);$$
\item[iv)] If $\overline{\partial}^4f=0,$ we have
$$u(z,\overline{z})=u_0(z)+\overline{z}f(z)-\frac{\overline{z}^2}{2}\overline{\partial}(f)+\frac{\overline{z}^3}{6}\overline{\partial}^2(f)-\frac{\overline{z}^4}{24}\overline{\partial}^3(f);$$
\end{itemize}
where $u_0$ belongs each time to some suitable space of holomorphic functions on $\Omega$ depending on the order of polyanalyticity $n=1,2,3,4$.
\end{example}

\begin{remark}
Let $f:\mathbb{C}\longrightarrow \mathbb{C}$ be such that $\overline{\partial}^n(f)=0$; the computations considered in the previous examples suggest a general expression for the particular solution of the  Problem \ref{pb-p} which will be of the form
\begin{equation}
\displaystyle u(z)=\sum_{k=1}^n \beta_k \overline{z}^k\overline{\partial}^{k-1}(f)(z),\quad \forall z\in \mathbb{C},
\end{equation}
where $\beta_k$ are suitable coefficients. We will prove that these coefficients are of the form
\begin{equation}
\beta_k=\frac{(-1)^{k+1}}{k!},\quad\forall k\geq 1.
\end{equation}

\end{remark}

Now we introduce some Sobolev type spaces of polyanalytic functions which will be useful for our purpose:
\begin{definition}[Hörmander-Sobolev type spaces] \label{HSTS}
Let $\Omega$ be an open set of $\mathbb{C}$,  and let $p$ be a subharmonic function in $\Omega$. We say that a measurable function $f:\Omega\longrightarrow \mathbb{C}$ belongs to the space $W^n_{p}(\Omega)$, $n=1,2,...,$ if for every $k=0,..., n-1$ we have
\begin{equation}
\displaystyle \int_{\Omega}|\overline{\partial}^kf(z)|^2e^{-p(z)}d\lambda(z)=M_k(f)<+\infty,
\end{equation}
with $d\lambda$ being the Lebesgue measure on $\Omega.$
We define also the Hörmander-Sobolev type spaces of polyanalytic functions of order $n=1,2,...$ by setting

\begin{equation}
\mathsf{A}^n_p(\Omega):=H_n(\Omega)\cap W^n_{p}(\Omega).
\end{equation}
\end{definition}

In the next result we present particular solutions for Problem \ref{pb-p}:
\begin{proposition}\label{PSestimate}
Let $n\in \mathbb{N}^*$ and $f\in \mathsf{A}^n_p(\Omega)$, i.e,  $f$ is a polyanalytic function of order $n$ such that  for every $ k=0,1,...,n-1$ we have

\begin{equation}
\displaystyle \int_{\Omega}|\overline{\partial}^{k} f(z)|^2e^{-p(z)}d\lambda(z)<+\infty.
\end{equation}  The function \begin{equation}\label{upsol}
u(z,\overline{z}):=-\displaystyle \sum_{k=1}^n \frac{(-1)^{k}}{k!}\overline{z}^k\overline{\partial}^{k-1}(f)(z), \quad \forall z\in \Omega.
\end{equation}
is a solution of $\overline{\partial}u=f.$
Moreover, the following estimate holds
\begin{equation}
\displaystyle \int_{\Omega}\frac{|u(z)|^2}{(1+|z|^2)^{2n}}e^{-p(z)}d\lambda(z) \leq n\sum_{k=1}^n \int_{\Omega}|\overline{\partial}^{k-1} f(z)|^2e^{-p(z)}d\lambda(z)<+\infty.
\end{equation}
\end{proposition}
\begin{proof}
We apply Leibniz rule with respect to $\overline{\partial}$-operator to the expression given by formula \eqref{upsol} and obtain
\begin{align*}
	\overline{\partial}u(z,\overline{z}) &=-\displaystyle \sum_{k=1}^n \frac{(-1)^{k}}{k!}\overline{\partial} \left(\overline{z}^k\overline{\partial}^{k-1}(f)\right)(z)\\
	&=-\displaystyle \sum_{k=1}^n \frac{(-1)^{k}}{k!} \left(k\overline{z}^{k-1}\overline{\partial}^{k-1}(f)+\overline{z}^k\overline{\partial}^{k}(f)\right) \\
	&=-\sum_{h=0}^{n-1}\frac{(-1)^{h+1}}{h!}\overline{z}^h\overline{\partial}^h(f)-\sum_{k=1}^{n}\frac{(-1)^k}{k!}\overline{z}^k\overline{\partial}^k(f)\\
	&=f(z)+ \sum_{k=1}^{n-1}\frac{(-1)^k}{k!}\overline{z}^k\overline{\partial}^k(f)-\sum_{k=1}^{n-1}\frac{(-1)^k}{k!}\overline{z}^k\overline{\partial}^k(f)\\
	&=f(z),
\end{align*}
which proves that $u$ is a particular solution of the $\overline{\partial}$-problem if the datum $f$ is polyanalytic of order $n$. As for the estimates, we recall the Cauchy-Schwarz inequality: for $a_1,...,a_n\geq 0$ we have

\begin{equation}
\displaystyle \left(\sum_{i=1}^na_i\right)^2\leq n\left(\sum_{i=1}^na_i^2\right).
\end{equation}
In our case this gives $$\displaystyle |u(z)| \leq \sum_{k=1}^n\frac{|z|^k}{k!}\left|\overline{\partial}^{k-1}(f)(z)\right| . $$
In particular we obtain

\begin{align*}
	\displaystyle |u(z)|^2 &\leq \left(\sum_{k=1}^n\frac{|z|^k}{k!}\left|\overline{\partial}^{k-1}(f)(z)\right| \right)^2 \\
	&\leq n \sum_{k=1}^n\frac{|z|^{2k}}{(k!)^2} \left|\overline{\partial}^{k-1}(f)(z)\right|^2.
\end{align*}

We observe that for every $k=1,...,n$ we have $$\displaystyle (1+|z|^2)^{2n}=\sum_{\ell=0}^{n}{2n \choose \ell}|z|^{2\ell}\geq |z|^{2k}.$$
Thus, it follows that
$$\displaystyle \frac{|z|^{2k}}{(1+|z|^2)^{2n}}\leq 1, \quad k=1,...,n.$$
Hence we obtain
\begin{align*}
	\displaystyle \int_{\Omega}\frac{|u(z)|^2}{(1+|z|^2)^{2n}}e^{-p(z)}d\lambda(z) &\leq n\sum_{k=1}^n\frac{1}{(k!)^2} \int_{\Omega} \frac{|z|^{2k}}{(1+|z|^2)^{2n}}|\overline{\partial}^{k-1} f(z)|^2e^{-p(z)}d\lambda(z)\\
	&\leq n \sum_{k=1}^n \frac{1}{(k!)^2}\int_{\Omega}|\overline{\partial}^{k-1} f(z)|^2e^{-p(z)}d\lambda(z)
	\\
	&\leq n \sum_{k=1}^n \int_{\Omega}|\overline{\partial}^{k-1} f(z)|^2e^{-p(z)}d\lambda(z).
\end{align*}

Finally, we have the following estimate
$$\displaystyle \int_{\Omega}\frac{|u(z)|^2}{(1+|z|^2)^{2n}}e^{-p(z)}d\lambda(z)\leq n\sum_{k=1}^n\frac{1}{(k!)^2} \int_{\Omega}|\overline{\partial}^{k-1} f(z)|^2e^{-p(z)}d\lambda(z)$$
$$\leq n \sum_{k=1}^n \int_{\Omega}|\overline{\partial}^{k-1} f(z)|^2e^{-p(z)}d\lambda(z).$$
\end{proof}

Now we will prove the converse of the previous result, thus obtaining a characterization for all solutions of the $\overline{\partial}$-problem for a given polyanalytic datum:
\begin{theorem}\label{GPADATUM}
Let $n\in \mathbb{N}^*$ and $f\in \mathsf{A}^n_p(\Omega)$, i.e,  $f$ is a polyanalytic function of order $n$ such that  for every $ k=0,1,...,n-1$ we have

\begin{equation}
\displaystyle \int_{\Omega}|\overline{\partial}^{k} f(z)|^2e^{-p(z)}d\lambda(z)<+\infty.
\end{equation}
 If $u$ is a solution of $\overline{\partial}u=f,$ satisfying the $L^2-$estimates,
then there exists a holomorphic function $u_0\in H(\Omega)$ such that \begin{equation}
u(z,\overline{z})=\displaystyle u_0(z)-\sum_{k=1}^{n}\frac{(-1)^{k}}{k!}\overline{z}^k\overline{\partial}^{k-1}(f)(z); \quad \forall z\in \Omega.
\end{equation}
Moreover, there exists a constant $C=C(n)>0$ depending on the order $n$ such that 

\begin{equation}
\displaystyle \int_{\Omega}\frac{|u_0(z)|^2}{(1+|z|^2)^{2n}}e^{-p(z)}d\lambda(z) \leq C(n) \left(\sum_{k=0}^{n-1} \int_{\Omega}|\overline{\partial}^{k} f(z)|^2e^{-p(z)}d\lambda(z)\right)<+\infty.
\end{equation}
\end{theorem}
\begin{proof}
We already presented in Example \ref{Expn4} the expressions of the solutions for $n=1,\ldots, 4$.
Now we prove a general formula for any $n\geq 1$. To this end, we use the fact that $f$ is a polyanalytic function of order $n$: by the poly-decomposition (\ref{polydecomposition}) there exist unique analytic functions $(f_k)_{0 \leq k\leq n-1}$ such that  \begin{equation}
\displaystyle f(z)=\sum_{k=0}^{n-1}\overline{z}^kf_k(z), \quad \forall z\in \Omega.
\end{equation}
Applying formula (13) of \cite{theodoresco1931derivee} we note that these analytic components $(f_k)_{0\leq k\leq n-1}$ can be expressed in terms of $f$ and of powers of $\overline{\partial}$-operator applied to $f$, so that we have:

$$\displaystyle f_k(z)=\frac{1}{k!}\sum_{s=0}^{n-k}\frac{(-1)^s}{s!}\overline{z}^s\overline{\partial}^{s+k}(f);\quad \forall k=1,..., n.$$

Now, since $u$ is a solution of the $\overline{\partial}$-problem associated to $f$, we know that $\overline{\partial}u(z,\overline{z})=f(z)$ and thus we get
$$\overline{\partial}^{n+1}u(z,\overline{z})=\overline{\partial}^{n}f(z)=0, \quad \forall z\in \Omega.$$
This shows that $u$ is a polyanalytic function of order $n+1$ which can be expressed in terms of $n$ analytic functions $(u_k)_{0 \leq k\leq n}$ uniquely determined, so that we can write
\begin{equation}\label{uform}
u(z,\overline{z})=u_0(z)+\overline{z}u_1(z)+\cdots+\overline{z}^n u_n(z)=u_0(z)+\sum_{k=1}^{n}\overline{z}^k u_k(z), \quad \forall z\in \Omega.
\end{equation}
Therefore, using the fact that $\overline{\partial}u=f$ we obtain
\begin{equation}
\displaystyle \sum_{k=1}^nk\overline{z}^{k-1}u_k(z)=\sum_{k=0}^{n-1}\overline{z}^{k}f_k(z), \quad \forall z\in \Omega.
\end{equation}
The uniqueness of the polyanalytic decomposition now allows us to identify the analytic components: 

 \begin{equation}
u_k(z)=\frac{1}{k}f_{k-1}(z)=\frac{1}{k!}\sum_{s=0}^{n-k}\frac{(-1)^s}{s!}\overline{z}^s\overline{\partial}^{s+k-1}(f),\quad \forall 1 \leq k\leq n.
\end{equation}
Inserting the previous expression of $(u_k(z))_{1\leq k \leq n}$ in formula \eqref{uform} we obtain

\begin{equation}
\displaystyle u(z)=u_0(z)+\sum_{k=1}^{n}\sum_{s=0}^{n-k}\frac{(-1)^s}{s!k!}\overline{z}^{s+k}\overline{\partial}^{s+k-1}(f)(z), \quad \forall z\in \Omega,
\end{equation}
which leads to
\begin{align*}
\displaystyle  u(z)&= u_0(z)+\sum_{\ell=1}^{n}\frac{(-1)^\ell}{\ell !}\overline{z}^\ell\overline{\partial}^{\ell-1}f(z) \left(1+\sum_{k=1}^{\ell} (-1)^k{\ell \choose k}-1\right)\\
&=u_0(z)+\sum_{\ell=1}^{n}\frac{(-1)^\ell}{\ell !}\overline{z}^\ell\overline{\partial}^{\ell-1}f(z) \left((1-1)^\ell-1\right)\\
&=u_0(z)-\sum_{\ell=1}^{n}\frac{(-1)^\ell}{\ell !}\overline{z}^\ell\overline{\partial}^{\ell-1}f(z).\\
\end{align*}

%$$\displaystyle u(z)=u_0(z)+\sum_{\ell=1}^{n}\frac{(-1)^\ell}{\ell !}\overline{z}\overline{\partial}^{\ell-1}f(z) \left(1+\sum_{k=1}^{\ell} (-1)^k{\ell \choose k}-1\right)$$

%Write $$u=u_0+\overline{z}u_1+...+\overline{z}^n u_n,$$
%we should use $$\overline{\partial}u=f,$$ and apply the unicity of the poly-decomposition combined with Teodoresco formulas ????

%We know by formula 13 in Todoresco thesis  that
%In particular, we observe also $$f_0(z)=f(z)-\sum_{k=1}^{n-1}\overline{z}^{k}f_k(z)$$
%Then, we use the fact that $$\overline{\partial}u=f,$$
%and use the unicity of poly-decomposition to get \begin{equation}
%u_k(z)=\frac{1}{k}f_{k-1}(z),\quad k\geq 1.
%\end{equation}
%Hence, \begin{equation}
%u(z)=\displaystyle u_0(z)+\sum_{k=1}^{n}\frac{\overline{z}^k}{k}f_{k-1}(z).
%\end{equation}
%Hence, we obtain
%$$u(z)=u_0(z)+\overline{z}f_0(z)+\sum_{k=2}^n \frac{\overline{z}^k}{k}f_{k-1}(z).$$
%Thus,

%$$u(z)=u_0(z)+\overline{z}f(z)-\sum_{k=2}^{n}\overline{z}^kf_{k-1}(z)+\sum_{k=2}^n\frac{\overline{z}^k}{k}f_{k-1}(z).$$
%Then, $$u(z)=u_0(z)+\overline{z}f(z)-\sum_{k=2}^{n}\overline{z}^k(\frac{k-1}{k})f_{k-1}(z).$$
%Now, we note that
%$$\frac{f_{k-1}(z)}{k}=\frac{1}{k!}\sum_{s=0}^{n-(k-1)}\frac{(-1)^s}{s!}\overline{z}^s\overline{\partial}^{s+k-1}(f);\quad k=2,..., n$$
%So, $$u(z)=u_0(z)+\overline{z}f(z)-\sum_{k=2}^{n}\overline{z}^k(\frac{k-1}{k!})\left(\sum_{s=0}^{n-(k-1)}\frac{(-1)^s}{s!}\overline{z}^s\overline{\partial}^{s+k-1}(f)\right).$$
%At this stage, we just need to replace $(f_{k-1})_{k\geq 2}$ by its expression and continue the computations.....\textbf{ (See Teoderesco fomula)}
Setting $\displaystyle \Psi_f(z):=\sum_{\ell=1}^{n}\frac{(-1)^\ell}{\ell !}\overline{z}^\ell\overline{\partial}^{\ell-1}f(z)$ and 
$$u_0(z)=u(z,\overline{z})+\Psi_f(z),\quad \forall  z\in \Omega,$$
we have
$$|u_0(z)|^2\leq2(|u(z,\overline{z})|^2+|\Psi_f(z))|^2),\quad \forall z\in \Omega.$$
Using the estimate obtained in Proposition \ref{PSestimate} we deduce that there exist some constants $C=C(n)>0$ such that

\begin{align*}
\displaystyle  \int_{\Omega}\frac{|u_0(z)|^2}{(1+|z|^2)^{2n}}e^{-p(z)}d\lambda(z) &\leq 2\left( \int_{\Omega}\frac{|u(z,\overline{z})|^2}{(1+|z|^2)^{2n}}e^{-p(z)}d\lambda(z) + \int_{\Omega}\frac{|\Psi_f(z)|^2}{(1+|z|^2)^{2n}}e^{-p(z)}d\lambda(z) \right) \\
&\leq 2n\left( \int_\Omega |f(z)|^2e^{-p(z)}d\lambda(z)+ \sum_{k=1}^n \int_{\Omega}|\overline{\partial}^{k-1} f(z)|^2e^{-p(z)}d\lambda(z)  \right) \\
&\leq C(n) \left(\sum_{k=0}^{n-1} \int_{\Omega}|\overline{\partial}^{k} f(z)|^2e^{-p(z)}d\lambda(z)\right).\\
\end{align*}

\end{proof}
We now turn our attention to another interesting class of polyanalytic function spaces, that arises when $p(z)=|z|^{\rho}$ for $\rho > 1.$
\begin{definition}
Let $\rho>1$ and $n\geq 1$. We define

$$\mathcal{A}^n_\rho(\mathbb{C}):=\left\lbrace{g \in H_n(\mathbb{C}) ;\quad \exists A,B>0:\quad |g(z)| \leq A \exp(B|z|^{\rho})}\right\rbrace, $$

%and

%$$\mathcal{A}^{n}_{\rho,0}(\mathbb{C}):=\left\lbrace{g \in H_n(\mathbb{C}) ;\quad \forall \varepsilon>0 \textbf{ } \exists A_\varepsilon>0:\quad |g(z)| \leq A_\varepsilon \exp(\varepsilon B|z|^{\rho})}\right\rbrace. $$
\end{definition}
\begin{proposition}
Let $n\geq1$ and $f_0,...,f_{n-1}\in \mathcal{A}_\rho^1(\mathbb{C})$. Then, the function defined by $$f(z)=\displaystyle \sum_{k=0}^{n-1}\overline{z}^kf_k(z),$$ belongs to the space $\mathcal{A}^{n}_{\rho}(\mathbb{C})$.
\end{proposition}
\begin{proof} It is clear that $f$ is polyanalytic of order $n$. Moreover,
we have
$$\displaystyle |f(z)|\leq \sum_{k=0}^{n-1}|z|^k|f_k(z)|\leq \sum_{k=0}^{n-1}|z|^kA_k\exp(B_k|z|^\rho)$$
Thus, by setting $A=\max(A_0,...,A_{n-1})$ and $B=\max(B_0,...,B_{n-1})$ we obtain $$|f(z)|\leq \left(\sum_{k=0}^{n-1}|z|^k \right)A\exp(B|z|^\rho) \leq A(\sum_{k=0}^{n-1} k!)\exp\left((B+1)|z|^\rho \right).$$
Hence, there exist two constants $A',B'>0$ which depend on the order of polyanalyticity $n$ so that we have
$$|f(z)| \leq A'\exp(B'|z|^p),\quad \forall z\in \mathbb{C}.$$
\end{proof}

\subsection{Examples: particular solutions for specific polyanalytic data}
In this final subsection we apply Theorem \ref{GPADATUM} in order to compute particular solutions of the $\overline{\partial}$-problem related to different examples of polyanalytic data. As a first example we prove the following
\begin{proposition}
Let $n\geq 1$ and consider the function defined by $f_n(z)=|z|^{2(n-1)}$ for every $z\in\Omega$. We have $f\in H_n(\Omega)$, and a particular solution of the problem $\overline{\partial}u=f$ is given by
$$\displaystyle u(z)=\frac{\overline{z}}{n}f_n(z).$$
\end{proposition}
\begin{proof}
We observe that $f_n$ is a polyanalytic function of order $n$. Then, by Theorem \ref{GPADATUM} it is clear that a particular solution of the $\overline{\partial}$-problem in this case is
 \begin{equation}
u(z):=-\displaystyle \sum_{k=1}^n \frac{(-1)^{k}}{k!}\overline{z}^k\overline{\partial}^{k-1}(f_n)(z); \quad \forall z\in \Omega.
\end{equation}
Since
$$\displaystyle \overline{\partial}^j(\overline{z}^p)=\frac{\Gamma(p+1)}{\Gamma(p-j+1)}\overline{z}^{p-j},\quad j=0,1,...$$
and

$$\displaystyle \sum_{k=1}^n (-1)^k{n \choose k}=-1,$$
we have $$\displaystyle \overline{\partial}^{k-1}(\overline{z}^{n-1})=\frac{\Gamma(n)}{\Gamma(n-k+1)}\overline{z}^{n-k},\quad k=1,2,...$$

Hence, $$\overline{\partial}^{k-1}(f_n)(z)=\frac{\Gamma(n)}{\Gamma(n-k+1)}z^{n-1}\overline{z}^{n-k},\quad k=1,2,... .$$
We finally have
$$u(z)=-\displaystyle \sum_{k=1}^n \frac{(-1)^{k}}{k!}\overline{z}^k\frac{\Gamma(n)}{\Gamma(n-k+1)}z^{n-1}\overline{z}^{n-k}=-\frac{1}{n}\overline{z}f_n(z) \sum_{k=1}^n (-1)^k{n \choose k}=\frac{\overline{z}}{n}f_n(z).$$
\end{proof}
The well-known complex Hermite polynomials introduced by Ito in \cite{Ito} form a very interesting class of polyanalytic functions. They are very useful, as they provide an orthogonal basis for the $L^2$ space on $\mathbb{C}$ with respect to the classical Gaussian measure since the following integral formula holds
$$\displaystyle \int_{\mathbb{C}}H_{m,n}(z,\overline{z})\overline{H_{m',n'}(z,\overline{z})}e^{-|z|^2}d\lambda(z)=\pi m!n! \delta_{((m,n);(m',n'))},$$
where $d\lambda(z)$ is the usual Lebesgue measure on $\mathbb{C}$.
\begin{definition}

Let $m,n=0,1, \cdots$, the complex Hermite polynomials are defined by:
\begin{equation}\label{Hmn}
H_{m,n}(z,\overline{z}):=\displaystyle \sum_{k=0}^{\min{(m,n)}}(-1)^k k! {m\choose k}{n \choose k}z^{m-k}\overline{z}^{n-k}.
\end{equation}
\end{definition}

\begin{remark}
In alternative, these polynomials can be defined via the Rodrigues formula:
\begin{equation}
H_{m,n}(z,\overline{z}):=(-1)^{m+n}e^{|z|^2}\frac{\partial^{m+n}}{\partial \overline{z}^m\partial z^n}(e^{-|z|^2}); \quad m,n=0,1,...
\end{equation}
\end{remark}
We have the following:
\begin{proposition}
Let $p,q=0,1,...$ and let $H_{p,q}(z,\overline{z})$ be the complex Hermite polynomials \eqref{Hmn}. Then, a particular solution of the $\overline{\partial}$-problem
$$\overline{\partial}u(z,\overline{z})=H_{p,q}(z,\overline{z}),$$

 is given by

$$u(z,\overline{z})=-\frac{1}{q+1}\sum_{k=1}^{q+1}(-1)^k{ q+1\choose k} \overline{z}^kH_{p, q+1-k}(z,\overline{z}),\quad \forall z\in \mathbb{C}.$$
\end{proposition}
\begin{proof}
We observe that $H_{p,q}$ is polyanalytic of order $q+1$. Moreover, Ito proved the following (see \cite{Ito}) 
$$\displaystyle \overline{\partial}H_{p,q}(z,\overline{z})=qH_{p,q-1}(z,\overline{z}),\quad q=1,2,...$$
Thus, we obtain the following

$$\displaystyle \overline{\partial}^j(H_{p,q}(z,\overline{z}))=\frac{\Gamma(q+1)}{\Gamma(q-j+1)}H_{p, q-j}(z,\overline{z}),\quad j=0,1,...$$
Then, in particular for $j=k-1$ we get

$$\displaystyle \overline{\partial}^{k-1}H_{p,q}(z)=\frac{\Gamma(q+1)}{\Gamma(q-(k-1)+1)}H_{p, q-(k-1)}(z,\overline{z})=\frac{\Gamma(q+1)}{\Gamma(q-k+2)}H_{p, q+1-k}(z,\overline{z}),\quad k=1,2,...$$
So, we obtain
$$\displaystyle \overline{\partial}^{k-1}H_{p,q}(z)=\frac{q!}{(q-k+1)!}H_{p, q+1-k}(z,\overline{z}),\quad k=1,2,...$$

Hence,
\begin{align*}
u(z,\overline{z})=& -\sum_{k=1}^{q+1}\frac{(-1)^k}{k!}\frac{q!}{(q-k+1)!}\overline{z}^kH_{p, q+1-k}(z,\overline{z})\\
= &-\sum_{k=1}^{q+1}\frac{(-1)^k}{(q+1)}{q+1\choose k} \overline{z}^kH_{p, q+1-k}(z,\overline{z}) \\
\end{align*}

%{
%\color{blue}
%Hence, we obtain

%$$u(z)=-\frac{1}{n}\left(\sum_{k=0}^{n} (-1)^k{ n\choose k} \overline{z}^kH_{p, n-k}(z,\overline{z}) -H_{p,n}(z,\overline{z}) \right).$$

%Finally, we note that

%$${ n\choose k} H_{p, n-k}(z,\overline{z})=\frac{1}{k!}\overline{\partial}^kH_{p,n}(z,\overline{z}).$$

%Thus, developing the computation we get

%$$u(z)=-\frac{1}{n}\sum_{k=1}^{n} \frac{(-\overline{z})^k\overline{\partial}^k}{k!} \left(H_{p,n}(z,\overline{z}) \right)$$
%We need to consider the following operator
%$$\mathcal{D}_n:=\sum_{k=1}^{n} \frac{(-\overline{z})^k\overline{\partial}^k}{k!}.$$
%Then, it holds that $$u=-\frac{1}{n}\mathcal{D}_{n}(f_{n}).$$

%We note that by Prop 3.1 (G) we have $$\displaystyle H_{p,q}(\sqrt{2}z,\sqrt{2}\overline{z})=\sqrt{2}^{p+q}q!\sum_{k=0}^q\frac{H_{0,k}(z,\overline{z})}{k!}\frac{H_{p,q-k}(z,\overline{z})}{(q-k)!}$$
%Moroever, we have the following  $$H_{0,k}(z,\overline{z})=\overline{z}^k.$$
%Then, we obtain
%$$u(z)=-\frac{1}{n}\left(\sum_{k=0}^{n} (-1)^k{ n\choose k} \overline{z}^kH_{p, n-k}(z,\overline{z}) -H_{p,n}(z,\overline{z}) \right)=.....$$
%}
\end{proof}

\begin{remark}
The particular solution $u(z,\overline{z})$ corresponding to the complex Hermite datum $H_{p,q}$ can be expressed as follows
$$u(z,\overline{z})=-\frac{1}{q+1}\mathcal{D}_{q+1}(H_{p,q+1});$$
where $$\mathcal{D}_{\ell}:=\sum_{k=1}^{\ell} \frac{(-1)^k}{k!}\overline{z}^k\overline{\partial}^k.$$
\end{remark}
We recall that $F_n(z,w)$ denote the reproducing kernel of the polyanalytic Fock space of order $n=1,2,...$, namely
\begin{equation}\label{Knn}
F_n(z,w)=F_{n,w}(z)=e^{z\overline{w}}\displaystyle \sum_{k=0}^{n-1}\frac{(-1)^k}{k!}{n \choose k+1}|z-w|^{2k}, \quad \forall z\in\mathbb{C}.
\end{equation}
Let us denote by $L^\alpha_n(x)$ the generalized Laguerre polynomials defined by
$$\displaystyle L^\alpha_m(x)=\displaystyle \sum_{k=0}^{m}\frac{(-x)^k}{k!}{ m+\alpha \choose m-k }. $$
The polyanalytic Fock kernel $F_n(z,w)$ given by \eqref{Knn} can be expressed in terms of the generalized Laguerre polynomials as follows
\begin{equation}\label{FpL}
F_n(z,w)=e^{z\overline{w}}L^1_{n-1}(|z-w|^2), \quad \forall z,w\in \mathbb{C}.
\end{equation}
\begin{lemma}\label{LF}
 Let $n\geq 1$ and $s=0,1,\cdots, n-1$ then we have
\begin{equation}
\displaystyle \overline{\partial}^{s}F_n(z,w)=(-1)^{s}(z-w)^{s}e^{z\overline{w}}\sum_{u=0}^{n-s-1}\frac{(-1)^{u}}{u!} {n\choose u+s+1}|z-w|^{2u} , \quad \forall z,w\in \mathbb{C}.
\end{equation}
%\begin{equation}
%\displaystyle \overline{\partial}^{k-1}F_n(z,w)=(-1)^{k+1}(z-w)^{k-1}F_{n-k%+1}(z,w), \quad k=1,2,\cdots, n.
%\end{equation}
In particular we have
\begin{equation}
\overline{\partial}^{n-1}F_n(z,w)=(-1)^{n-1}(z-w)^{n-1}F_1(z,w), \quad \forall z,w\in \mathbb{C}.
\end{equation}
\end{lemma}
\begin{proof}
Let $0 \leq s\leq n-1$, starting from the definition of the kernel $F_n(z,w)$ it is clear that
\begin{align*}
\displaystyle \overline{\partial}^sF_n(z,w)&=e^{z\overline{w}}\displaystyle \sum_{k=0}^{n-1}\frac{(-1)^k}{k!}{n \choose k+1}\overline{\partial}^s|z-w|^{2k}\\
&=e^{z\overline{w}}\displaystyle \sum_{k=0}^{n-1}\frac{(-1)^k}{k!}{n \choose k+1}(z-w)^k\overline{\partial}^s(\overline{z}-\overline{w})^{k}\\
\end{align*}
At this stage we use the fact that $\overline{\partial}^s(\overline{z}-\overline{w})^{k}=0$ if $k<s$ and $\overline{\partial}^s(\overline{z}-\overline{w})^{k}=\frac{k!}{(k-s)!}(\overline{z}-\overline{w})^{k-s}$ if $k\geq s$. In particular, this leads to

\begin{align*}
\displaystyle \overline{\partial}^sF_n(z,w)&=e^{z\overline{w}}\displaystyle \sum_{k=s}^{n-1}\frac{(-1)^k}{k!}{n \choose k+1} (z-w)^k  \frac{k!}{(k-s)!}(\overline{z}-\overline{w})^{k-s}   \\
&=e^{z\overline{w}}\displaystyle \sum_{k=s}^{n-1}\frac{(-1)^k}{(k-s)!}{n \choose k+1} (z-w)^k (\overline{z-w})^{k-s}\\
\end{align*}
We use the change of index $u=k-s$ and get
$$\displaystyle \overline{\partial}^sF_n(z,w)=(-1)^s(z-w)^se^{z\overline{w}}\sum_{u=0}^{n-s-1}\frac{(-1)^{u}}{u!} {n\choose u+s+1}|z-w|^{2u}.$$
Finally, by taking $s=n-1$ in the previous calculations we obtain
$$\displaystyle \overline{\partial}^{n-1}F_n(z,w)=(-1)^{n-1}(z-w)^{n-1}e^{z\overline{w}}=(-1)^{n-1}(z-w)^{n-1}F_1(z,w).$$
\end{proof}
\begin{corollary}
For every $n\geq 1$ we have
\begin{equation}
\displaystyle \overline{\partial_{z}}^{n-1} \partial_{w}^{n-1}F_n(z,w)=F_1(z,w), \quad \forall z,w\in \mathbb{C}.
\end{equation}
\end{corollary}
\begin{proof}
We already know by Lemma \ref{LF} that
$$\displaystyle \overline{\partial_z}^{n-1}F_n(z,w)=(-1)^{n-1}(z-w)^{n-1}F_1(z,w),$$
thus if we act with the operator $\partial_{w}^{n-1}$ we get

\begin{align*}
\displaystyle \partial_{w}^{n-1}\overline{\partial_z}^{n-1}F_n(z,w)&=(-1)^{n-1}\partial_{w}^{n-1}\left((z-w)^{n-1}\right) F_1(z,w)\\
&=F_1(z,w).\\
\end{align*}
So, taking the mixed derivatives up to the order $n-1$ with respect to both the variables $w$ and $\overline{z}$ applied to the polyanalytic Fock kernel $F_n(z,w)$ we re-obtain the classical Fock kernel $F_1(z,w)$ and this ends the proof.
\end{proof}

We observe that using our approach a particular solution of the $\overline{\partial}$-problem associated to the poly-Fock kernel datum $F_n$ can be expressed in terms of the generalized Laguerre polynomials as follows:
\begin{theorem}\label{pp}
For every fixed parameter $w\in \mathbb{C}$ we set  $f(z)=F_n(z,w)$. Then, a particular solution of the $\overline{\partial}$-problem $\overline{\partial}u=f$ is given by the following formula
\begin{equation}
 \displaystyle u_n(z,\overline{z})=e^{z\overline{w}}\sum_{s=1}^{n} \frac{\overline{z}^s}{s!}(z-w)^{s-1}L^{s}_{n-s}(|z-w|^2).
\end{equation}
\end{theorem}
\begin{proof}
We know by Proposition \ref{PSestimate} that a particular solution to the problem is given by
$$u_n(z,\overline{z})=\displaystyle- \sum_{s=1}^n\frac{(-1)^s}{s!}\overline{z}^s\overline{\partial}^{s-1}(f)(z), \quad z\in\mathbb{C}$$
So, we use Lemma \ref{LF} to compute $\overline{\partial}^{s-1}(f)(z)$ for $s\geq 1$ and get

\begin{align*}
\displaystyle \overline{\partial}^{s-1}(f)(z)&=\overline{\partial}^{s-1}F_n(z,w)\\
&=(-1)^{s-1}(z-w)^{s-1}e^{z\overline{w}}\sum_{u=0}^{n-s}\frac{(-1)^{u}}{u!}{n \choose u+s} |z-w|^{2u}.\\
\end{align*}
Thus, taking $\alpha=s$ and $m=n-s$ we observe the following relation with the generalized Laguerre polynomials
$$\displaystyle \overline{\partial}^{s-1}(f)(z)=(-1)^{s-1}(z-w)^{s-1}e^{z\overline{w}}L_{n-s}^{s}(|z-w|^2),\quad \forall s=1,\cdots, n.$$
 Therefore, we continue the computations and deduce that the particular solution $u$ is given by
 $$\displaystyle u_n(z,\overline{z})=e^{z\overline{w}}\sum_{s=1}^{n} \frac{\overline{z}^s}{s!}(z-w)^{s-1}L^{s}_{n-s}(|z-w|^2).$$

%$$\displaystyle \overline{\partial}^{s-1}(f)(z)=\overline{\partial}^{s-1}F_n(z,w)=(-1)^{s-1}(z-w)^{s-1}\sum_{u=0}^{n-s}\frac{(-1)^{u}}{u!}{n \choose u+s} |z-w|^{2u}$$
%$$ \displaystyle u_n(z,\overline{z})=e^{z\overline{w}}\sum_{s=1}^{n} \frac{\overline{z}^s}{s!}(z-w)^{s-1}\left(\sum_{u=0}^{n-s}\frac{(-1)^{u}}{u!}{n \choose u+s}|z-w|^{2u}\right).$$
%{\color{red}

%we note that $$\displaystyle \overline{\partial}^\ell|z-w|^{2k}=\frac{k!}{(k-\ell)!}(z-w)^\ell|z-w|^{2(k-\ell)}; \quad \ell=0,1,...$$

%Hence for $\ell=s-1$ we get
%$$\displaystyle \overline{\partial}^{s-1}|z-w|^{2k}=\frac{k!}{(k-s+1)!}(z-w)^{s-1}|z-w|^{2(k-s+1)}; \quad s=1,2,...$$

%Thus, we have
%\begin{equation}
%\overline{\partial}^{s-1}(f)(z)=e^{z\overline{w}}(z-w)^{s-1}\displaystyle \sum_{k=0}^{n-1}\frac{(-1)^k}{k!}{n \choose k+1} \frac{k!}{(k-s+1)!}|z-w|^{2(k-s+1)};\quad s=1,2,...
%\end{equation}
%}
\end{proof}

\begin{corollary}
With the above notations, for every $n \geq 2$ we have
\begin{equation}
u_n(z,\overline{z})=\overline{z}F_n(z,w)+F_1(z,w)\sum_{s=2}^{n}\frac{\overline{z}^s}{s!}(z-w)^{s-1}L_{n-s}^{s}(|z-w|^2).
\end{equation}
\end{corollary}
\begin{proof}
This is a direct application of Theorem \ref{pp} combined with formula \eqref{FpL}.
\end{proof}
\section*{Acknowledgments}
Daniel Alpay thanks the Foster G. and Mary McGaw Professorship in Mathematical Sciences, which supported this research. Daniele C. Struppa is grateful to the Donald Bren Presidential Chair in Mathematics. Kamal Diki thanks the Grand Challenges Initiative (GCI) at Chapman University, which supported this research.

%\bibliographystyle{plain}
%\bibliographystyle{plain}
%\bibliography{Hormanderpoly}

\begin{thebibliography}{10}
\bibitem{abramowitz1964handbook}
M. Abramowitz and I.A. Stegun.
\newblock  Handbook of mathematical functions with formulas, graphs, and
  mathematical tables, volume~55.
\newblock  US Government printing office, 1964.

\bibitem{abreu2010sampling}
L.D. Abreu.
\newblock Sampling and interpolation in Bargmann-Fock spaces of polyanalytic
  functions.
\newblock {\em Applied and Computational Harmonic Analysis}, {\bf 29} (3) (2010):287--302.

\bibitem{abreu2012super}
L.D. Abreu.
\newblock Super-wavelets versus poly-Bergman spaces.
\newblock {\em Integral Equations and Operator Theory}, {\bf 73} (2) (2012):177--193.

\bibitem{abreu2014function}
L.D. Abreu and H.G. Feichtinger.
\newblock Function spaces of polyanalytic functions.
\newblock In {\em Harmonic and complex analysis and its applications}, pages
  1--38. Springer, 2014.

\bibitem{agranovsky2011characterization}
M.L. Agranovsky.
\newblock Characterization of polyanalytic functions by meromorphic extensions
  from chains of circles.
\newblock {\em Journal d'Analyse Math{\'e}matique}, {\bf 113} (1) (2011) :305--329.

\bibitem{AhernBruna} P. Ahern, and J. Bruna.
\newblock Maximal and area integral characterizations of Hardy-Sobolev spaces in the unit ball of $\mathbb C^ n$.
\newblock {\it Revista matemática iberoamericana} {\bf 4} (1) (1988): 123-153.


\bibitem{alpay2020generalized}
D. Alpay, P. Cerejeiras, and U. Kaehler.
\newblock Generalized Fock space and moments.
\newblock {\em arXiv preprint arXiv:2005.08085}, 2020.

\bibitem{alpay2022generalized}
D. Alpay, P. Cerejeiras, and U. K{\"a}hler.
\newblock Generalized white noise analysis and topological algebras.
\newblock {\em Stochastics}, pages 1--33, 2021.

\bibitem{HFspace} D. Alpay, F. Colombo, K. Diki, I. Sabadini, and D.C. Struppa.
\newblock A Hörmander-Fock type space.
\newblock Preprint 2022.

\bibitem{askour2000explicit}
N. Askour, A. Intissar, and Z. Mouayn.
\newblock Explicit formulas for reproducing kernels of generalized Bargmann
  spaces on $\mathbb{C}^n$.
\newblock {\em Journal of Mathematical Physics}, {\bf 41} (5) (2000):3057--3067.

\bibitem{Balk1991}
M.~Balk.
\newblock {\em Polyanalytic functions}.
\newblock Akademie-Verlag, Berlin., 1991.

\bibitem{balk1997polyanalytic}
M.~Balk.
\newblock Polyanalytic functions and their generalizations.
\newblock In {\em Complex Analysis I}, pages 195--253. Springer, 1997.

\bibitem{begehr2002orthogonal}
H. Begehr.
\newblock Orthogonal decompositions of the function space $L_2 (\overline{D}; \mathbb{C})$.
\newblock  J. Reine Angew. Math. {\bf 549} (2002): 191--219.

\bibitem{berenstein1993complex}
C.A Berenstein and D.C. Struppa.
\newblock Complex analysis and convolution equations.
\newblock In {\em Several Complex Variables V}, pages 1--108. Springer, 1993.

\bibitem{berenstein1979new}
C.A Berenstein and B.A. Taylor.
\newblock A new look at interpolation theory for entire functions of one
  variable, {\it Adv. Math.} {\bf 33} (2) (1979): 109-143.


\bibitem{berndtsson2010introduction}
B. Berndtsson.
\newblock An introduction to things $\overline{\partial}$.
\newblock {\em Analytic and Algebraic Geometry, McNeal}, pages 7--76, 2010.

\bibitem{el2014primer}
O. El-Fallah, K. Kellay, J. Mashreghi, and T. Ransford.
\newblock {\em A primer on the Dirichlet space}, volume 203.
\newblock Cambridge University Press, 2014.

\bibitem{haslinger2014d}
F. Haslinger.
\newblock The d-bar Neumann problem and Schr{\"o}dinger operators.
\newblock In {\em The d-bar Neumann Problem and Schr{\"o}dinger Operators}. De
  Gruyter, 2014.

\bibitem{hedenmalm2015hormander}
H. Hedenmalm.
\newblock On H{\"o}rmander’s solution of the $\overline{\partial}$-equation, I, 
\newblock {\em Mathematische Zeitschrift}, {\bf 281} (1) (2015): 349--355.

\bibitem{hormander1966l2}
L. H{\"o}rmander.
\newblock $L^2$ estimates and existence theorems for the $\overline{\partial}$ operator.
\newblock {\em Acta Math}, {\bf 113} (1966): 89--152.

\bibitem{hormander1967generators}
L. H{\"o}rmander.
\newblock Generators for some rings of analytic functions.
\newblock {\em Bulletin of the American Mathematical Society}, {\bf 73} (6) (1967): 943--949.

\bibitem{hormander2007notions}
L. H{\"o}rmander.
\newblock {\em Notions of convexity}.
\newblock Springer Science \& Business Media, 2007.

\bibitem{Ito}
K. It{\^o}.
\newblock Complex multiple Wiener integral.
\newblock In {\em Japanese Journal of Mathematics: Transactions and Abstracts},
  volume~22, pages 63--86. The Mathematical Society of Japan, 1952.

\bibitem{kolossoff1908problemes}
G.V. Kolossoff.
\newblock Sur les probl{\`e}mes d’{\'e}lasticit{\'e} {\`a} deux dimensions.
\newblock {\em CR Acad. Sci}, {\bf 148} (1908):1242--1244.

\bibitem{theodoresco1931derivee}
N. Th{\'e}odoresco.
\newblock La d{\'e}riv{\'e}e ar{\'e}olaire et ses applications {\`a} la
  physique math{\'e}matique.
\newblock 1931.

\bibitem{vasilevski1999structure}
N.L. Vasilevski.
\newblock On the structure of Bergman and poly-Bergman spaces.
\newblock {\em Integral Equations and Operator Theory}, {\bf 33} (4) (1999): 471--488.


\end{thebibliography}

%\begin{thebibliography}{10}

%\end{thebibliography}

\end{document}